\documentclass[12pt]{article}
\usepackage[left=1in,top=1in,right=1in,bottom=1in]{geometry}

\usepackage[english]{babel}
\usepackage{amsmath,amsthm,amsfonts,amssymb,epsfig,color,authblk, hyperref}
\usepackage{bbm,booktabs,float,mathtools,siunitx,pgfplots,tikz}

\newtheorem{thm}{Theorem}[section]

\newtheorem{conj}{Conjecture}

\newtheorem*{ass}{Assumption}

\newcommand{\E}{\mathbf{E}}

\newcommand{\mG}{\mathcal{G}}
\newcommand{\PP}{\mathbf{P}}
\newcommand{\wPP}{\widetilde{\mathbf{P}}}
\newcommand{\R}{\mathbb{R}}

\newcommand{\ii}{{\rm i}}
\newcommand{\oo}{{\rm o}}
\newcommand{\rd}{{\rm d}}

\newcommand{\re}{{\rm Re}}
\newcommand{\im}{{\rm Im}}
\def\one{\ensuremath\mathbbm{1}}

\title{\Large Evidence of Random Matrix Corrections for the Large Deviations of Selberg's Central Limit Theorem}

\author[1]{\normalsize  E. Amzallag}
\author[2,3]{\normalsize L.-P. Arguin}
\author[3,4]{\normalsize E. Bailey}
\author[2]{\normalsize K. Hui}
\author[2]{\normalsize R. Rao}
\affil[1]{\footnotesize \it Department of Mathematics, City College of New York, CUNY, New York, NY}
\affil[2]{\footnotesize  \it Department of Mathematics, Baruch College, CUNY, New York, NY}
\affil[3]{\footnotesize  \it Department of Mathematics, CUNY Graduate Center, New York, NY}
\affil[4]{\footnotesize  \it Department of Mathematics, University of Bristol, UK}

\date{April 16, 2021}

\begin{document}

\maketitle
\begin{abstract}
  Selberg's central limit theorem states that the values of $\log|\zeta(1/2+\ii \tau)|$, where $\tau$ is a uniform random variable on $[T,2T]$, is distributed like a Gaussian random variable of mean $0$ and standard deviation $\sqrt{\frac{1}{2}\log \log T}$. It was conjectured by Radziwi{\l}{\l} that this breaks down for values of order $\log\log T$, where a multiplicative correction $C_k$ would be present at level $k\log\log T$, $k>0$. This constant should be equal to the leading asymptotic for the $2k^{th}$ moment of $\zeta$, as first conjectured by Keating and Snaith using random matrix theory. In this paper, we provide numerical and theoretical evidence for this conjecture. We propose that this correction has a significant effect on the distribution of the maximum of $\log|\zeta|$ in intervals of size $(\log T)^\theta$, $\theta>0$. The precision of the prediction enables the numerical detection of $C_k$ even for low $T$'s of order $T=10^8$. A similar correction appears in the large deviations of the Keating-Snaith central limit theorem for the logarithm of the characteristic polynomial of a random unitary matrix, as first proved by F\'eray, M\'eliot and Nikeghbali. 
\end{abstract}

\section{Introduction and Main Results}
\subsection{Introduction}
The large values of the Riemann zeta function $\zeta: \mathbb C\to \mathbb C$ on the critical line $\re (s)=1/2$ play an important role in number theory. 
There are several conjectures describing its purported behavior.
For example, the {\it moment conjecture} gives the precise asymptotics of the moments of the function on the interval $[T,2T]$ (see for example~\cite{keasna00a, ivic85, titchmarsh86}):
\begin{conj}[Moment Conjecture]
  \label{conj: moments}
  For $k\geq 0$, as $T\rightarrow\infty$,
  \begin{equation}
    \label{eqn: moments}
    \frac{1}{T}\int_0^T|\zeta(1/2+\ii t)|^{2k}\rd t\sim C_{k}(\log T)^{k^2}.
  \end{equation}
\end{conj}
Lower bounds agreeing with the predicted leading exponent are known unconditionally (see \cite{radsou13, heasou20} and the earlier works of \cite{heabro81, ram94}). Consistent upper bounds on the level of the leading exponent are known unconditionally for $0\leq k\leq 2$ (from the work of \cite{hearadsou19}), and for all $k\geq0$ conditionally on the Riemann hypothesis \cite{har13b, sou09}. This paper is predominantly about the constants $C_k$.

At the level of the constant, Conjecture~\ref{conj: moments} is only proved in the case $k=1$, by Hardy and Littlewood with $C_1=1$, and the case $k=2$ by Ingham with $C_2=\frac{1}{(2\pi)^2}$ \cite{harlit18, ing26}. The constants $C_k$, $k> 0$, have been conjectured in \cite{keasna00a} using random matrix theory to be of the following form
\begin{equation}
  \label{eqn: C}
  C_k=a_k\cdot f_k,
\end{equation}
where
\begin{align}
  a_k&= \prod_{p \text{ primes}}\left(1-\frac{1}{p}\right)^{k^2}\sum_{m=0}^\infty\left(\frac{\Gamma(k+m)}{m!\Gamma(k)}\right)^2p^{-m},\label{coeff_arith}
\end{align}
and
\begin{align}
  f_k&=\frac{\mG^2(1+k)}{\mG(1+2k)}.\label{coeff_rmt}
\end{align}
Here, $\mathcal{G}$ denotes the Barnes $\mathcal{G}-$function. An alternative approach using Dirichlet series yields the same conjecture \cite{diagolhof03}.

As pointed out in \cite{rad11}, the constants $C_k$ should also appear in the large deviations of Selberg's central limit theorem. 
The theorem asserts that for $\sigma^2_T=\frac{1}{2}\log\log T$
\begin{equation}
  \label{eqn: selberg}
  \PP\Big(\log |\zeta(1/2+\ii\tau)|>\sigma_T\cdot V\Big)\sim \int_V^\infty\frac{e^{-\frac{x^2}{2}}}{\sqrt{2\pi}}\rd  x, \quad V\in \R.
\end{equation}
However, for $V$ of the order of the variance, the proposed correction is:
\begin{conj}[Radziwi\l\l's Conjecture\footnote{There is a small typo in the statement of the conjecture in the original paper where $V\sim k\sqrt{\log\log T}$. A factor of $2$ in the square root is missing.}]
  \label{conj: maks}
  If $\tau$ is uniformly distributed on $[T,2T]$ and $\sigma^2_T=\frac{1}{2}\log\log T$, then for $V\sim k\sqrt{2\log\log T}$ we have for $k>0$,
  \begin{equation}
    \label{eqn: maks conj}
    \PP\Big(\log |\zeta(1/2+\ii\tau)|>\sigma_T\cdot V\Big)\sim 
    C_{k}\int_V^\infty\frac{e^{-\frac{x^2}{2}}}{\sqrt{2\pi}}\rd  x.
  \end{equation}
\end{conj}
The conjecture is plausible since if one biases the choice of $\tau$ by the value $|\zeta(1/2+\ii \tau)|^{2k}$, the left-hand side of \eqref{eqn: maks conj} becomes:
\begin{equation}
  \PP\Big(\log |\zeta(1/2+\ii\tau)|>\sigma_T\cdot V\Big)=\E[|\zeta(1/2+\ii\tau)|^{2k}]\cdot \wPP\Big(\log |\zeta(1/2+\ii\tau)|>\sigma_T\cdot V\Big),
\end{equation}
where $\wPP$ is defined by $\frac{\rd \wPP}{\rd \PP}=\frac{|\zeta(1/2+\ii \tau)|^{2k}}{\E[|\zeta(1/2+\ii \tau)|^{2k}]}$.
The distribution of $\log |\zeta(1/2+\ii\tau)|$ has been recently proved to be Gaussian under $\wPP$ as in Selberg's theorem, and with the same variance \cite{faz21}.
The constant $C_k$ appears naturally from the asymptotics \eqref{eqn: moments} of the moments.

\subsection{Results}
The main objective of this paper is to provide more evidence that the moment correction $C_k$, predicted by random matrix theory, should be present in \eqref{eqn: maks conj}.
On the theoretical side, it was proved by F\'eray, M\'eliot, and Nikeghbali that a similar correction naturally appears in the large deviations of the Keating-Snaith central limit theorem for the characteristic polynomial of the circular unitary ensemble (CUE):
\begin{thm}[Theorem 7.5.1 \cite{FerMelNik2016}]
  \label{thm: LD KS}
  Let $P_N(\theta)=\det(I-e^{\ii \theta}U)$ be the characteristic polynomial of an $N\times N$ random matrix $U$ sampled under the Haar measure $\mathbb P_{\mathcal{U}(N)}$ on the unitary group $\mathcal{U}(N)$.
  Write $Q_2(N)$ for the second cumulant of $\log|P_N(\theta)|$. 
  Then we have for $V\sim k \sqrt{2\log N}$, $k\geq 0$, and any $\theta\in [0,2\pi)$
    \[
    \mathbb P_{\mathcal{U}(N)}\Big(\log |P_N(\theta)|>\sqrt{Q_2(N)}\cdot V\Big)\sim 
    f_k \int_V^\infty\frac{e^{-\frac{x^2}{2}}}{\sqrt{2\pi}}\rd x, \quad \text{ where } f_k=\dfrac{\mG^2(1+k)}{\mG(1+2k)}.
    \]
\end{thm}
The statement of Theorem~\ref{thm: LD KS} is precisely the random matrix analogue of~\eqref{eqn: maks conj}.  Using the usual dictionary (see, e.g., \cite{keasna00a, katsar99b, cfkrs05}), one compares the unitary characteristic polynomial $P_N(\theta)$ with $\zeta(1/2+\ii t)$, and by comparing densities of eigenvalues and zeros of zeta, the matrix size $N$ corresponds to a height $\log(t/2\pi)$.

Large deviations in various ranges were also considered by Hughes et al.~\cite{hugkeaoco01}.  In that paper, they show that $\log|P_N(\theta)|/A(N)$ satisfies a large deviation principle, for various ranges of $A(N)$. The appropriate range in the context of Theorem~\ref{thm: LD KS} is \emph{moderate} deviations: $\sqrt{\log{N}}\ll A(N)\ll N$, and it is shown in~\cite{hugkeaoco01} that the rate function is either quadratic or linear, depending on the precise growth of $A(N)$. 
Similar results have also been proved  in the general context of $\beta$-ensembles and Wigner matrices in \cite{DorEic2013,DorEic2013Cumu}.
Theorem \ref{thm: LD KS} differs in that it examines a very particular form of $A(N)$ and derives the resulting precise constant multiple of the Gaussian. 
The proof of Theorem \ref{thm: LD KS} was done in \cite{FerMelNik2016} in the general context for mod-$\phi$ convergence.
For completeness, we provide the detailed computation following the result of Keating \& Snaith in Section \ref{sect: proof}.

The statement of Theorem~\ref{thm: LD KS} equally holds for the imaginary part of $\log P_N(\theta)$, with the appropriate change to the leading order correction\footnote{In this case, one would instead find that $f_k=\mathcal{G}(1-k)\mathcal{G}(1+k)$, the leading coefficient of the $2k^{th}$ moment of $e^{\ii\im\log P_N(\theta)}$.}.  Additionally, the statement of Theorem \ref{thm: LD KS} can be generalized  to other classical compact groups, where the distribution should again be Gaussian but with a different correction $g_k$ (corresponding to the relevant matrix group moment), see for example~\cite{FerMelNik2016}. From such calculations one could deduce conjectures akin to Conjecture~\ref{conj: maks} for symplectic and orthogonal families of $L$-functions (cf.~\cite{cfkrs05}).

Beyond the theoretical evidence from Theorem~\ref{thm: LD KS}, we also provide numerical evidence of the presence of the correction $C_k$. 
We choose to investigate its effect on the maximum of the real part of the logarithm of the zeta function in short intervals, instead of directly testing \eqref{eqn: maks conj}.
The reason for this is that the presence of the correction leads to a very precise refinement of the Fyodorov-Hiary-Keating conjecture for the maximum of $\log |\zeta(1/2+\ii \tau+\ii h)|$ for $h$ in a short interval, \cite{fyohiakea12,fyokea14}. The conjecture was originally stated for {\it mesoscopic} intervals, that is, intervals of size $2\pi(\log T)^{\theta}$ for $\theta\leq 0$.
The statement was adapted in \cite{argouirad19} to {\it macroscopic} intervals of size $2\pi(\log T)^\theta$, for $\theta>0$, but fell short of capturing the order one terms.
It turns out that the correction $C_k$ has a measurable effect on the recentering of the maximum. More precisely, we propose the following refinement:
\begin{conj}
  \label{conj: IID}
  Consider a fixed $\theta>0$.
  If $\tau$ is uniformly distributed on $[T,2T]$, then we have
  \begin{equation}
    \label{eqn: to test}
    \max_{|h|\leq \pi(\log T)^{\theta}} \log |\zeta(1/2+\ii(\tau+h))|
    =\sqrt{1+\theta} \log \log T-\frac{1}{4\sqrt{1+\theta}}\log\log\log T+ \mathcal G_{\theta,T},
  \end{equation}
  where  $(\mathcal G_{\theta,T}, T\geq 1)$ is a sequence of random variables converging in distribution to a Gumbel random variable $\mathcal G_\theta$ with $\PP(\mathcal G_\theta\leq x)=\exp(-e^{-\frac{1}{\beta}(x-m)})$ and parameters
  \begin{equation}
    \begin{aligned}
      \label{eqn: m beta}
      \beta&=\beta(\theta)=\frac{1}{2\sqrt{1+\theta}}\\
      m&=m(\theta)=(0.06537\dots)+\beta^2\log C_{\sqrt{1+\theta}}-\frac{\beta^2}{2}\Big(\log (1+\theta)- \log (4\pi)\Big). 
    \end{aligned}
  \end{equation}
\end{conj}
As can be seen from the leading order of \eqref{eqn: to test}, the relevant regime of large deviation at a given $\theta$ is $\sqrt{1+\theta}\log\log T$.
Together with Conjecture \ref{conj: maks}, this leads naturally to the choice $k=\sqrt{1+\theta}$ for $C_k$ in \eqref{eqn: m beta}.  
The precise numerical constant appearing in the definition of $m$ in \eqref{eqn: m beta} is the Meissel-Mertens constant divided by $4$, see Equation \eqref{eqn: variance B}. 

The upshot of Conjecture \ref{conj: IID} is a very precise prediction to order one for the maximum of $\log |\zeta(1/2+\ii(\tau+h))|$, including very good control of the finite-size effects, that can be compared to the numerical data. The high precision of the conjecture to order one is the saving grace here, as the factors $\log \log T$, $\log\log\log T$, and $\mathcal G_{\theta, T}$ in \eqref{eqn: to test} remain essentially of the same order for all testable $T$'s (around $T=10^{23}$ seems to be the current computational limit). 
In particular, this spares us some of the difficulty of testing the moment conjecture, see \cite{hiaodl12}.

We chose to test Equations \eqref{eqn: to test} and \eqref{eqn: m beta} at $T=10^7, 10^8$ and $10^9$, where it is not costly to amass a good sample size for many $\theta$'s for $0\leq \theta \leq 3$. A snapshot of the results for the empirical mean of the maximum are given in Figure \ref{fig: table} and Figures \ref{fig: mean 7}, \ref{fig: mean 8} and \ref{fig: mean 9}. The main conclusion there is that the correction $C_k$ is necessary to fit the data. Details on the numerical experiments are given in Section \ref{sect: num}.

The case $\theta=0$ is special as Conjecture \ref{conj: IID} is not expected to hold. It was proposed in \cite{fyohiakea12,fyokea14} that the subleading order should instead be $-\frac{3}{4}\log\log T$. One then would expect the empirical mean to lie \emph{lower} than the prediction \eqref{eqn: to test}. 
In addition, the fluctuations should not be exactly Gumbel but a randomly shifted Gumbel. The effect of the random shift is such that the right tail of the distribution of the recentered maximum is not exponential, as for a pure Gumbel, but should be heavier: of the form $ye^{-\beta y}$. This would in effect \emph{increase} the contribution of the fluctuations to the mean. 
The two above corrections seem hard to observe numerically. One problem is that they are competing effects, which may mutually cancel. Secondly, there is the systematic problem that the standard deviation of the maximum is fairly large at $\theta=0$, as can be seen for example in Figure \ref{fig: mean 9}.
Theoretical progress to settle the Fyodorov-Hiary-Keating conjecture has recently been made in \cite{naj18,abbrs19, Har19,argbourad20}.
A continuous smoothing of the subleading order between $-\frac{1}{4\sqrt{1+\theta}}\log\log T$ to $-\frac{3}{4}\log\log T$ as $\theta\downarrow 0$ has been proposed in \cite{argdubhar21} by taking $\theta\sim (\log\log T)^{-\alpha}$, $0<\alpha<1$. This gives a subleading order of $-\frac{(1+2\alpha)}{4}\log\log T$.
Again, this interpolation seems hard to capture numerically as the standard deviation of the maximum is large for small $\theta$.

The paper is structured as follows.
We provide a proof of Theorem \ref{thm: LD KS} in Section \ref{sect: proof}, following the work of Keating and Snaith~\cite{keasna00a}. 
The details on how the theoretical prediction based on Conjecture \ref{conj: IID} is generated are given in Section \ref{sect: theoretical prediction}.
The conjecture is derived in Section \ref{sect: derivation} using basic extremal value theory, assuming Conjecture \ref{conj: maks} and reasonable properties of $\zeta$.
We comment on the control of the finite-size effects and on the numerical computations of the $C_k$'s for all $0\leq k\leq 2$ in Section \ref{sect: finite-size}. Numerical experiments are discussed in Section \ref{sect: num}.
\\

\noindent{\bf Acknowledgements.}
We thank M. Radziwi{\l\l} for valuable comments on the first version of the paper. We are also grateful to P. Bourgade and A. Nikeghbali for pointing out to us that Theorem \ref{thm: LD KS} first appeared in \cite{FerMelNik2016}.
L.-P. A. gratefully acknowledges the support from the grant NSF CAREER~DMS-1653602. K. H. and R. R. were also financially supported in part by this grant.  E. B. thanks the Heilbronn Institute for Mathematical Research for support.


\section{Proof of Theorem \ref{thm: LD KS}}
\label{sect: proof}
For completeness, we provide the proof of Theorem \ref{thm: LD KS} following Keating \& Snaith \cite{keasna00a}.
Recall that we write $P_N(\theta)=\det(I-e^{\ii \theta}U)$ for the characteristic polynomial of a random $N\times N$ unitary matrix.  Let $V\sim k\sqrt{2\log N}$.  We will show that
\begin{equation}
  \mathbb P_{\mathcal{U}(N)}\left(\frac{\log|P_N(\theta)|}{\sqrt{Q_2(N)}}> V\right)\sim f_k \int_V^\infty \frac{e^{-\frac{x^2}{2}}}{\sqrt{2\pi}}\rd x,
\end{equation}
where
\begin{equation}\label{eq:rmt_coeff}
  f_k=\frac{\mathcal{G}^2(1+k)}{\mathcal{G}(1+2k)},
\end{equation}
and $Q_2(N)$ is the second cumulant of $\log|P_N(\theta)|$.  Keating and Snaith calculated the asymptotic form of the cumulants $Q_j(N)$,
\begin{align}
  Q_1(N)&=0\\
  Q_2(N)&=\frac{1}{2}\log N+\frac{1}{2}(\gamma+1)+\frac{1}{24 N^2}-\frac{1}{80 N^4}+\mathcal{O}\left(\frac{1}{N^6}\right)\label{eq:second_cumulant}\\
  Q_m(N)&=(-1)^m\frac{2^{m-1}-1}{2^{m-1}}\left(\Gamma(m)\zeta(m-1)-\frac{(m-3)!}{N^{m-2}}\right)+\mathcal{O}\left(\frac{1}{N^{m-1}}\right),
\end{align}
for $m\geq 3$. 

Let the probability density function of $\log|P_N(\theta)|$ be
\begin{equation}
  \rho_N(x)= \mathbb{E}_{\mathcal{U}(N)}[\one\{\log|P_N(\theta)|=x\}],
\end{equation}
where the average is taken over the unitary group with respect to Haar measure. Additionally define the moment generating function for $\log|P_N(\theta)|$ to be
\begin{equation}\label{eq:mgf}
  M_N(s)=\sum_{j=0}^\infty\frac{\mathbb{E}_{\mathcal{U}(N)}\left[(\log|P_N(\theta)|)^j\right]}{j!}s^j.
\end{equation}
Using the Selberg integral, Keating and Snaith~\cite{keasna00a} determined a finite $N$ formula for $M_N(s)$, valid for all real $\theta$ and $\re(s)>-1$
\begin{equation}\label{eq:mgf2}
  M_N(s)=\mathbb{E}_{\mathcal{U}(N)}[|P_N(\theta)|^s]=\prod_{j=1}^N\frac{\Gamma(j)\Gamma(j+s)}{\Gamma^2(j+s/2)}.
\end{equation}
Due to the rotational invariance of the Haar measure on $\mathcal{U}(N)$, the right hand side of \eqref{eq:mgf2} is independent of $\theta$.  From \eqref{eq:mgf2}, one deduces that, as $N\rightarrow\infty$, $M_N(s)\sim f_k N^{k^2}$, where $f_k$ is given by \eqref{eq:rmt_coeff}. 

Define
\begin{equation}\label{eq:twisted_pdf}
  \tilde{\rho}_N(x)= \sqrt{Q_2(N)}\rho_N(\sqrt{Q_2(N)}x),
\end{equation}
the probability density function of $\log|P_N(\theta)|$ rescaled.
Keating and Snaith determined (cf.~\cite{keasna00a}, Equation (53)) that
\begin{equation}\label{eq:rhotilde}
  \tilde{\rho}_N(x)=\frac{1}{\sqrt{2\pi}}e^{-\frac{x^2}{2}}\left(1+\sum_{m=3}^\infty \frac{A_m(N)}{Q_2(N)^{\frac{m}{2}}}\sum_{p=0}^m\binom{m}{p}x^p\mathcal{E}(m,p)\right),
\end{equation}
where
\begin{equation}
  \label{eq:parity sum}
  \mathcal{E}(m,p)
  =
  \begin{cases}
    \ii^{m-p}(m-p-1)!!,&m-p\text{ even},\\
    0,&m-p\text{ odd}
  \end{cases},
\end{equation}
and where the terms $A_m(N)$ are determined combinatorially from the cumulants and Equation \eqref{eq:twisted_pdf}. For example,
\begin{align*}
  A_3(N)&=\frac{Q_3(N)}{3!}\\
  A_4(N)&=\frac{Q_4(N)}{4!}\\
  A_5(N)&=\frac{Q_5(N)}{5!}\\
  A_6(N)&=\frac{Q_6(N)}{6!}+\frac{1}{2!}\frac{Q_3(N)^2}{(3!)^2}.
\end{align*}
As $N\rightarrow\infty$, $A_m(N)$ approaches a constant, see~\cite{keasna00a}.
Evaluating Equation \eqref{eq:rhotilde} at $x=\tfrac{k\log N}{\sqrt{Q_2(N)}}$ then gives
\begin{equation}
  \tilde{\rho}_N\left(\frac{k\log N}{\sqrt{Q_2(N)}}\right)
  =\frac{1}{\sqrt{2\pi}}e^{-\frac{1}{2}\frac{k^2\log^2 N}{Q_2(N)}}
  \left(1+\sum_{m=3}^\infty \frac{A_m(N)}{Q_2(N)^{\frac{m}{2}}}\sum_{p=0}^m\binom{m}{p}\left(\frac{k\log N}{\sqrt{Q_2(N)}}\right)^{p}
  \mathcal E(m,p)\right).
\end{equation}
Note that $x\sim k\sqrt{2\log N}$ as $N\rightarrow\infty$.
Also, for large $N$, the sum over $p$ is dominated at $p=m$, hence
\begin{align}
  \tilde{\rho}_N\left(\frac{k\log N}{\sqrt{Q_2(N)}}\right)&\sim\frac{1}{\sqrt{2\pi}}e^{-\frac{1}{2}\frac{k^2\log^2 N}{Q_2(N)}}\left(1+\sum_{m=3}^\infty A_m(N)\left(\frac{k\log N}{Q_2(N)}\right)^{m}\right)\\
  &\sim\frac{1}{\sqrt{2\pi}}e^{-k^2\log N}e^{k^2(\gamma+1)}\left(1+\sum_{m=3}^\infty A_m(N)(2k)^m\right),\label{eq:pdf_limit}
\end{align}
using Equation \eqref{eq:second_cumulant}.  

Finally, recall from the comment following Equation \eqref{eq:mgf2} that
\begin{equation}
  f_k= \lim_{N\rightarrow\infty}\frac{M_N(2k)}{N^{k^2}}.
\end{equation}
Then, since the cumulants $Q_m(N)$ are related to the Taylor coefficients of $\log M_N(s)$ (cf. Equation~\eqref{eq:mgf}),
\begin{equation}
  M_N(s)=\exp\left(\sum_{m=1}^\infty\frac{Q_m(N)}{m!}s^m\right),
\end{equation}
we have
\begin{align}
  \frac{M_N(2k)}{N^{k^2}}
  &=\exp\left(2k^2Q_2(N)-k^2\log N+\sum_{m= 3}^\infty\frac{ Q_m(N)}{m!}(2k)^m\right)\\
  &\sim\exp\left(k^2(\gamma+1)+\sum_{m=3}^\infty\frac{ Q_m(N)}{m!}(2k)^m\right)\\
  &=\exp\left(k^2(\gamma+1)\right)\left(1+\sum_{m=3}^\infty A_m(N)(2k)^m\right).\label{eq:alt_fk}
\end{align}
The final line follows by expanding the the exponential with the infinite sum in the exponent in the Taylor series, and then grouping terms according to power of $2k$.  The fact that this results in the weights $A_m(N)$ follows immediately from the combinatorial definition, see the discussion following \eqref{eq:parity sum} and~\cite{keasna00a}.  Hence $f_k$ is given by Equation \eqref{eq:alt_fk}, and the result follows by substituting the expression into \eqref{eq:pdf_limit}.

\section{Derivation of Prediction}
\label{sect: theoretical prediction}

\subsection{Derivation of Conjecture \ref{conj: IID}}
\label{sect: derivation}
In this section, we derive Conjecture \ref{conj: IID} based on the following assumptions:
\begin{ass}
  For $\tau$ a uniform random variable on $[T,2T]$, the stochastic process
  \[( |\zeta(1/2+\ii(\tau+h))|, |h|\leq \pi (\log T)^\theta)\]
  satisfies the following:
  \begin{enumerate}
  \item {\it Discretization}: the maximum over the interval $[-\pi(\log T)^{\theta}, \pi(\log T)^{\theta}]$ can be reduced to the maximum over a discrete set $\mathcal H_{\theta,T}$ corresponding to the midpoints between the zeros of $|\zeta(1/2+\ii(\tau+h))|$ on the interval.\label{ass1}
  \item {\it Independence}: The variables $|\zeta(1/2+\ii(\tau+h))|$, $h\in \mathcal H_{\theta, T}$, are independent. \label{ass2}
  \item {\it Large deviations of Selberg's central limit theorem}: Equation \eqref{eqn: maks conj} holds with $k=\sqrt{1+\theta}$.\label{ass3} 
  \end{enumerate}
\end{ass}
Assumption~\ref{ass1} is reasonable as the maximum should be achieved between two zeros. This can be rigorously established, see \cite{argouirad19}.
Assumption~\ref{ass2} cannot be exact, but it is likely a very good approximation. Indeed, the correlation between the values at $h$ and $h'$ decays very fast with the distance: as $|h-h'|^{-1}$, \cite{argouirad19}. Assumption~\ref{ass3} is the one to be tested. The choice of $k$ comes from the expected leading order of the maximum being $\sqrt{1+\theta}\log\log T$. 

From these assumptions, the derivation of the distribution of the maximum is a standard computation in extreme value theory. However, we shall need very good control of the finite-size effects to compare with numerics, so we include the details. The finite-size effects are discussed in the next section.

The number of zeros $\mathcal N(t)$ on $[0,t]$ is known to a very good level of precision thanks to the Riemann-von Mangoldt formula, see for example \cite{tsang84,ivic85},
\begin{equation}
  \mathcal N(t)=\frac{t}{2\pi}\log \frac{t}{2\pi e} + \frac{1}{\pi}{\rm Im} \log \zeta(1/2+\ii t) +\mathcal O(1).
\end{equation}
This implies that the number of zeros in the interval $[\tau-\pi(\log T)^{\theta},\tau+\pi(\log T)^{\theta}]$ is
\begin{equation}
  \label{eqn: N}
  \begin{aligned}
    N_{\theta,T}&=\mathcal N(\tau+\pi (\log T)^{\theta}) - \mathcal N(\tau-\pi (\log T)^{\theta})\\
    &=(\log T)^\theta \log \frac{\tau}{2\pi e}+\mathcal O((\log T)^\theta)\\
    &=(\log T)^{1+\theta}+\mathcal O((\log T)^\theta).
  \end{aligned}
\end{equation}
This will be the approximation for the cardinality of the discrete set $\mathcal H_{\theta,T}$.
The above implies
\begin{equation}
  \log N_{\theta,T}=(1+\theta)\log\log T+\mathcal O((\log T)^{-1}).
\end{equation}

Focussing on Assumption~\ref{ass3}, the right tail of the distribution of $\log |\zeta(1/2+\ii \tau)|$ is expected to be Gaussian with multiplicative correction $C_{k}$ and variance
\begin{equation}
  \label{eqn: sigma1}
  \sigma^2_T=\frac{1}{2}\sum_{p\leq T}\frac{1}{p}.
\end{equation}
In \cite{rad11}, Conjecture \ref{conj: maks} is stated for $\sigma_T^2=\frac{1}{2} \log\log T$. Whereas it is true that $\frac{1}{2}\sum_{p\leq T}\frac{1}{p}=(1+\oo(1))\log\log T$, the $\oo(1)$-term is needed for precise numerics. Mertens' second theorem asserts that, see for example \cite{Van17},
\begin{equation}
\label{eqn: variance B}
  \sum_{p\leq T}\frac{1}{p}=\log\log T+B+\mathcal O\left(\frac{e^{-\sqrt{0.175\log T}}}{(\log T)^{3/4}}\right),
\end{equation}
where $B=0.26149\dots$ is the Meissel-Mertens contant. Hence, from \eqref{eqn: sigma1}, the standard deviation is asymptotically 
\begin{equation}
  \label{eqn: sigma}
  \sigma_T=\sqrt{\tfrac{1}{2}\log\log T}+\frac{B}{2\sqrt{2\log\log T}}+ \mathcal O(e^{-\sqrt{0.1\log T}}).
\end{equation}

We are now ready to derive Conjecture \ref{conj: IID}. We use the shorthand notation $N=N_{\theta,T}$, $\sigma=\sigma_T$ and $C=C_{\sqrt{1+\theta}}$ (the moment coefficient).
Under Assumptions \ref{ass1}, \ref{ass2}, and \ref{ass3}, we have for any $Y>0$ with $\sigma Y\sim \sqrt{1+\theta}\log\log T$,
\begin{equation}
  \label{eqn: max1}
  \PP\left(\max_{|h|\leq \pi(\log T)^\theta}\log|\zeta(1/2+\ii(\tau+h))|>\sigma \ Y\right)
  =\left(1-\frac{NC\int_{Y}^{\infty} \frac{e^{-\frac{x^2}{2}}}{\sqrt{2\pi}}\rd x}{N}\right)^{N}.
\end{equation}
First, note that
\begin{equation}
  \label{eqn: E1}
  \left(1-\frac{a}{n}\right)^{n}=e^{-a}\cdot\left(1-\mathcal O\left(\frac{a^2}{n}\right)\right),
\end{equation}
so that
\begin{equation}
  \label{eqn: max2}
  \PP\left(\max_{|h|\leq \pi(\log T)^\theta}\log|\zeta(1/2+\ii(\tau+h))|>\sigma \ Y\right)\sim \exp\left(- N C\int_{Y}^{\infty} \frac{e^{-\frac{x^2}{2}}}{\sqrt{2\pi}}\rd x\right).
\end{equation}
The correct level of the maximum is obtained by choosing $Y$ for which the numerator in Equation \eqref{eqn: max1} is of order one. With this mind, consider $Y^\star$ the solution to the equation
\begin{equation}
  \label{eqn: Y}
  N C\int_{Y^{\star}}^{\infty} \frac{e^{-\frac{x^2}{2}}}{\sqrt{2\pi}}\rd x=1.
\end{equation}
Write $\sigma Y=\sigma \ Y^\star+y$, $y\in\R$.
With this notation, Equation \eqref{eqn: max2} becomes
\begin{equation}
  \PP\left(\max_{|h|\leq \pi(\log T)^\theta}\log|\zeta(1/2+\ii(\tau+h))|>\sigma \ Y^\star +y\right)\sim \exp\left(- G(y)\right),
\end{equation}
where 
\begin{equation}
  \label{eqn: G}
  G(y)=\dfrac{\int_{Y^\star +y/\sigma}^{\infty} e^{-\frac{x^2}{2}}\rd x}{\int_{ Y^\star}^{\infty} e^{-\frac{x^2}{2}}\rd x}.
\end{equation}

The quantity $\sigma Y^\star$ is the deterministic shift in Equation \eqref{eqn: m beta}. To see this, recall that standard Gaussian estimates give the asymptotics
\begin{equation}
  \label{eqn: estimate}
  \left(1-\frac{1}{Y^2}\right)\frac{1}{Y}\frac{e^{-\frac{Y^2}{2}}}{\sqrt{2\pi}}\leq\int_{Y}^{\infty} \frac{e^{-\frac{y^2}{2}}}{\sqrt{2\pi}}\rd y\leq\frac{1}{Y}\frac{e^{-\frac{Y^2}{2}}}{\sqrt{2\pi}}, \quad Y\geq 1.
\end{equation}
In other words, it is possible to write
\begin{equation}
  \label{eqn: E2}
  \int_{Y}^{\infty} \frac{e^{-\frac{y^2}{2}}}{\sqrt{2\pi}}\rd y=\left(1-\mathcal O\left(\frac{1}{Y^2}\right)\right)  \frac{1}{Y}\frac{e^{-\frac{Y^2}{2}}}{\sqrt{2\pi}} .
\end{equation}
Combining \eqref{eqn: Y} and \eqref{eqn: E2} gives the following equation for $Y^\star$:
\begin{equation}
  \label{eqn: Y2}
  \frac{ NC}{\sqrt{2\pi}}\frac{e^{-\frac{(Y^\star)^2}{2}}}{Y^\star}=1+\mathcal{O}\left(\frac{1}{(Y^\star)^2}\right).
\end{equation}
The solution can be approximated recursively. 
A first approximation omitting $\frac{1}{Y^\star}=e^{-\log Y^\star}$ yields $Y^\star\approx \sqrt{2\log N}$.
Writing $Y^\star=\sqrt{2\log N}+\delta$ in \eqref{eqn: Y2} gives an equation for $\delta$:
\begin{equation}
  \frac{\delta^2}{2}+\delta \sqrt{2\log N}+\frac{1}{2}\log\log N+\frac{1}{2}\log 2+\log \left(1+\frac{\delta}{\sqrt{2\log N}}\right)=\log C-\frac{1}{2}\log 2\pi.
\end{equation}
It is straightforward to solve this by expanding the logarithm to get
\begin{align}
  Y^\star
  &=\sqrt{2\log N}+\delta\\
  &=\sqrt{2\log N}-\frac{1}{2\sqrt{2\log N}}\Big(\log\log N+\log{4\pi}-2\log C\Big)+\mathcal O\left(\frac{\log\log N}{\log N}\right).\label{eqn: Y star}
\end{align}
It remains to note that Equation \eqref{eqn: N} and Equation \eqref{eqn: sigma} imply 
\begin{align}
  \sigma\sqrt{2\log N}&=\sqrt{1+\theta}\ \log\log T +\sqrt{1+\theta}\ \frac{B}{2}+\oo(1),\\
  \frac{\sigma}{2\sqrt{2\log N}}&=\frac{1}{4\sqrt{1+\theta}}\left(1+\mathcal O((\log\log T)^{-1})\right).
\end{align}
After multiplication by $\sigma$, the first two terms in Equation \eqref{eqn: Y star} give the leading correction for the maximum of $\log |\zeta|$:
$$
\sqrt{1+\theta}\log\log T-\frac{1}{4\sqrt{1+\theta}}\log\log\log T.
$$
The remaining terms with $\log C$, $B$, $\log (1+\theta)$ and $\log 4\pi$ amount to the deterministic shift $m$ in Equation \eqref{eqn: m beta}.

It remains to study the fluctuations around $Y^\star$. 
This is done by computing the asymptotics of the function $G$, defined by \eqref{eqn: G}.
Using the Gaussian estimate \eqref{eqn: estimate} again yields
\begin{equation}
  \label{eqn: G2}
  G(y)\sim e^{-\frac{Y^\star}{\sigma} y}\cdot \frac{e^{-\frac{y^2}{(2\sigma^2)}}}{1+\frac{y}{\sigma Y^\star}}.
\end{equation}
Equation \eqref{eqn: m beta} follows by noting that $e^{-y^2/(2\sigma^2)}\sim 1$ by \eqref{eqn: sigma}, and $1+\frac{y}{\sigma Y^\star}\sim 1$ for any fixed $y$ by \eqref{eqn: Y star}.
Moreover, Equations \eqref{eqn: N}, \eqref{eqn: sigma} and \eqref{eqn: Y star} imply
\begin{equation}
  \label{eqn: beta}
  \frac{Y^\star}{\sigma}\sim 2\sqrt{1+\theta}.
\end{equation}
Putting this back in \eqref{eqn: G2} gives the parameter $\beta$ of the Gumbel distribution in Equation \eqref{eqn: m beta}. 

\subsection{Discussion of the Finite-Size Effects}
\label{sect: finite-size}
In this section, we explain how the numerical predictions of Conjecture \ref{conj: IID} are obtained, including the computations of the coefficients $C_k$.

According to Conjecture \ref{conj: IID}, the mean of $\max_{|h|\leq \pi(\log T)^\theta}\log|\zeta(1/2+\ii(\tau+h))|$ consists of two terms: the deterministic recentering
\begin{equation}
  \sqrt{1+\theta}\log\log T-\frac{1}{4\sqrt{1+\theta}}\log\log \log T,
\end{equation}
and the expectation of the Gumbel random variable $\mathcal G_{\theta, T}$, which in the limit is
\begin{equation}
  \label{eqn: mean gumbel}
  \E[\mathcal G_\theta]=m+\beta \gamma, 
\end{equation}
where $\gamma$ is the Euler constant $\gamma=0.577\dots$, and with the parameters $m$ and $\beta$ given in Equation \eqref{eqn: m beta}.

To give an idea of the orders of magnitude in the problem considered, note that for $T=10^8$, we have
\begin{equation}
\label{eqn: T8}
\log\log 10^8=2.78\dots, \quad \log\log\log 10^8=1.07\dots\ .
\end{equation}
Scaling the simulations up to to $T=10^{23}$ wouldn't result in much additional precision, since there
\begin{equation}
\log\log 10^{23}=3.97\dots, \quad \log\log\log 10^{23}=1.39\dots \ .
\end{equation}
This modest gain would come at a substantial computational time cost. 
Despite the curse of iterated logarithms, it is possible to derive accurate numerical predictions thanks to a precise control of the lower order terms and their finite-size corrections. 

The computation of the deterministic shift $\sigma Y^\star$ follows the treatment in the last section. 
For the number of zeros $N_{\theta, T}$, we rely on Equation \eqref{eqn: N} and take
\begin{equation}
N_{\theta,T}\approx(\log T)^\theta \log \frac{T}{2\pi e}.
\end{equation}
The error term in \eqref{eqn: N} is of order $(\log T)^{-1}$, which is comparatively small. 
For the variance $\sigma_T^2$, note that the constant $B=0.26149\dots$ is fairly close to the value of $\log\log T$ in view of \eqref{eqn: T8}. 
For this reason, we use Equation \ref{eqn: sigma} in the computation. One might expect a quadratic correction to $\sum_{p\leq T}\frac{1}{p}$ of the form $\sum_{p\leq T}\frac{1}{8p^2}\approx \frac{1}{32}$, due to the expansion of the Euler product. However, this is negligible for our purpose.
The approximation of the product by the exponential in Equation   \ref{eqn: E1} also comes at a low cost, since the multiplicative error $(1-\mathcal O(a^2/n))$ is evaluated at $a=1$,  by design in \eqref{eqn: Y}, and at $n=N_{\theta,T}\approx (\log T)^\theta\geq 18.4$. This multiplicative error becomes of small shift in the exponential, allowing us to discard it.
The Gaussian estimate \eqref{eqn: Y2} is not quite precise enough for these fine numerics. Indeed, this leads to Equation \eqref{eqn: Y star}, with an error $\mathcal O\left(\frac{\log\log N_{\theta,T}}{\log N_{\theta,T}}\right)$. 
This error remains substantial at any $T$'s that are computationally reachable. For this reason, we approximate $Y^\star$ directly by numerically solving Equation \eqref{eqn: Y}. 
This takes care of all finite-size corrections to $\sigma Y^\star$. 

It remains to evaluate the mean of the Gumbel random variable $\mathcal G_{\theta, T}$.
As can be seen from Equation \eqref{eqn: G2}, the function $G$ converges to $e^{-2\sqrt{1+\theta}y}$, albeit very slowly. The term $e^{-y^2/2\sigma^2}$ is particularly problematic since $\sigma$ is of the order $\sqrt{\log\log T}$. 
To take care of this, instead of using $\beta \gamma$ in the mean of a Gumbel (cf. Equation~\eqref{eqn: mean gumbel}), we simply evaluate the mean by working directly with the function $G(y)$ in Equation \eqref{eqn: G}. 
The mean of the recentered random variable $\max_{|h|\leq \pi(\log T)^\theta}\log|\zeta(1/2+\ii(\tau+h))|-\sigma Y^\star$ can then be evaluated  using the cumulative distribution function $1-\exp(-G(y))$. In Equation \eqref{eqn: beta}, we also use $\frac{\sigma^2}{\sigma Y^\star}$ for $\beta$ instead of the limiting value $(2\sqrt{1+\theta})^{-1}$. For $\theta=3$, these considerations yield a mean of $0.17\dots$ compared to $0.14\dots$ for the limiting function $e^{-2\sqrt{1+\theta}y}$. In the same way, the variance of $\max_{|h|\leq \pi(\log T)^\theta}\log|\zeta(1/2+\ii(\tau+h))|$ in the limit should be the one of a Gumbel with parameter $\beta$, i.e.,
\begin{equation}
\label{eqn: std gumbel}
\text{Var}(\mathcal G_\theta)=\frac{\beta^2\pi^2}{6}.
\end{equation}
The finite-size effects of computing the variance with the CDF $1-\exp(-G(y))$ are tiny compared to using \eqref{eqn: std gumbel}, so we use Equation \eqref{eqn: std gumbel} for simplicity.

Finally, we turn to computing the moment coefficient $C_k$ (which recall appears in the definition of $m$, see~\eqref{eqn: m beta}).  For our purpose, we are interested in the value of $C_k=a_kf_k$ given by Equation \eqref{eqn: C} for $k\in(1,2)$, since $k=\sqrt{1+\theta}$ and for our numerics we take $\theta\in(0,3)$. Certain evaluations of $a_k$ and $f_k$ appear in the literature: see for example~\cite{hiaodl12}, where they compute $C_k$ for the first few integers $k$:
\begin{table}[H]
  \centering
  \begin{tabular}{c|c}
    \toprule
    $k$ & $C_k=a_kf_k$\\
    \midrule
    $1$ & $1$ \\
    $2$ & $\frac{1}{2\pi^2}\approx \num{5.066e-2}$   \\
    $3$ & $\num{5.708e-6}$\\
    $4$ & $\num{2.465e-13}$\\
    \bottomrule
  \end{tabular}
  \caption{Values of the leading order coefficient in the moment conjecture~\eqref{eqn: moments}. The exact values for $k=1,2$ are due to Hardy and Littlewood, and Ingham respectively~\cite{harlit18, ing26}.  The (truncated) numerical values for higher $k$ can be found in~\cite{hiaodl12}.} 
\end{table}
The numerical values of $C_k=a_kf_k$ presented in Section~\ref{sec:numerics} were computed\footnote{Computations were completed in SageMath \cite{sagemath}, version 9.1, using Python 3.7.} by evaluating \eqref{coeff_arith} and \eqref{coeff_rmt}. Figure~\ref{fig:coeff_plot} plots $C_k$ for $\theta\in(0,3)$ (i.e. $k\in(1,2)$).


  \caption{Graph of the conjectured leading order coefficient $C_{k}$ of the $2k^{th}$ moment of $|\zeta(1/2+it)|$, for $k=\sqrt{1+\theta}$. The horizontal dashed line is at $(2\pi^2)a_2f_2=1$.}
  \label{fig:coeff_plot}
\end{figure}

\section{Numerical Experiments}\label{sec:numerics}
\label{sect: num}
\subsection{About the Experiments}
Thanks to the precision of Conjecture \ref{conj: IID}, the prediction can be tested using fairly rudimentary numerical experiments.
The datasets were generated using Python 3.8. We employ SageMath's $lcalc$ function \cite{sagemath}, which uses Michael Rubinstein's $L$-function calculator.
 We also implement multiprocessing to expedite the run-time for large simulations.
 Samples were constructed at $T=10^7, 10^8$ and $10^9$ over the interval $0\leq \theta \leq 3$ for each $0.1$ increment for $T=10^7, 10^8$ and for each 0.25 increment for $T=10^9$.
 
For each height, $T$, the method consisted on generating $\mathcal S$ evaluations of $\tau$, for $\tau$ uniformly distributed on $[T,2T]$.  Hence, the interval $[\tau-\pi(\log T)^{\theta},\tau+\pi(\log T)^{\theta}]$ was discretized at every $\frac{2\pi}{\log(T/2\pi)}$, and the maximum over this discrete set of points computed.
The sample sizes $\mathcal S$ at each $\theta$ were $500$, $400$, $300$ for $T=10^7$, $10^8$ and $10^9$ respectively. 

\subsection{Results}

\begin{figure}[H]
\begin{center}
\includegraphics[height=8cm]{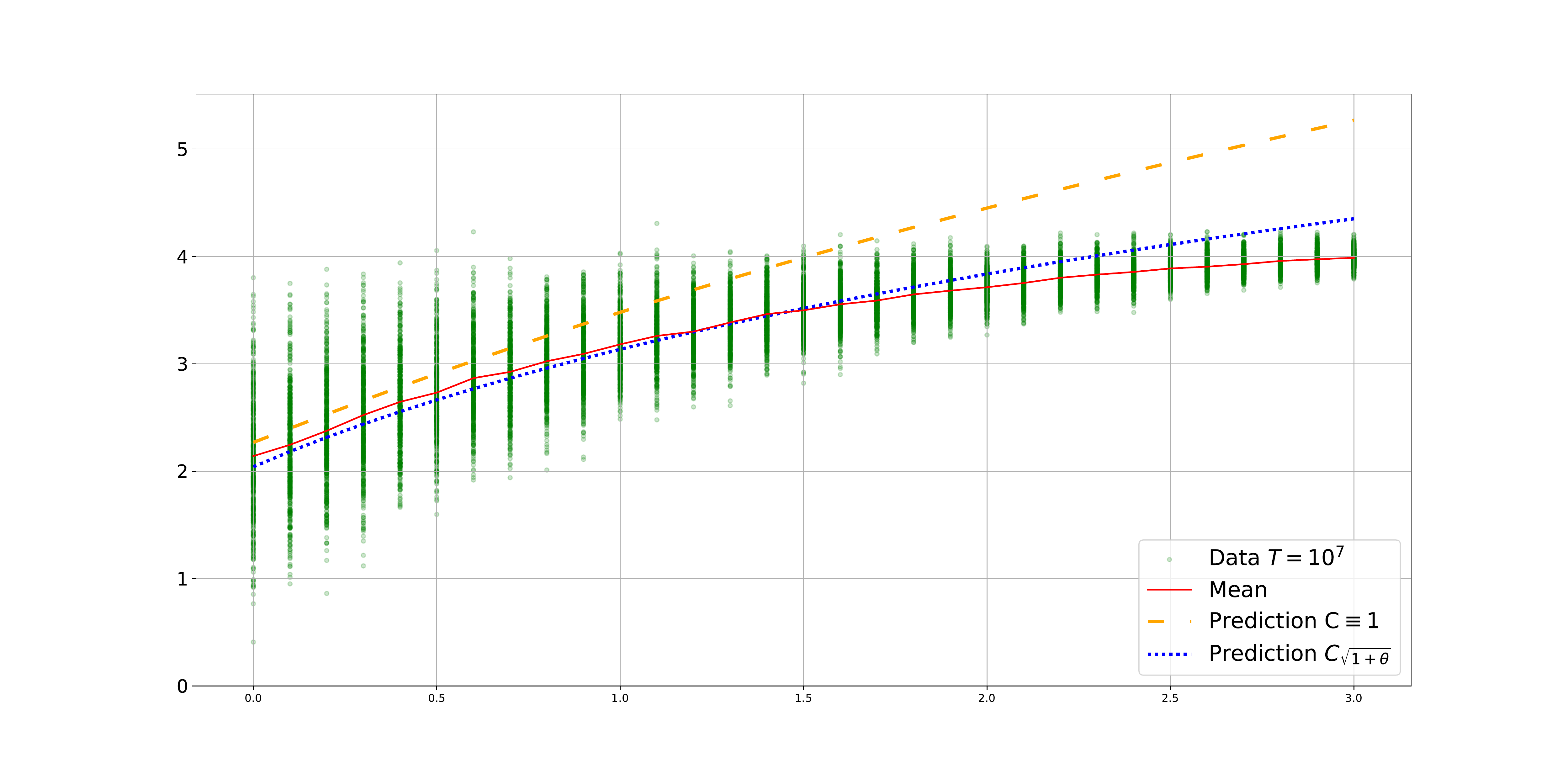}
\end{center}
\caption{The empirical  mean of the samples of $\max_{|h|\leq\pi (\log T)^\theta}\log |\zeta(1/2+\ii(\tau+h))|$ as a function of $\theta$, $0\leq \theta\leq 3$, with step size $0.1$ at $T=10^7$. The dotted lines correspond to the theoretical predictions for $C\equiv 1$ and for $C_{\sqrt{1+\theta}}$.}
\label{fig: mean 7}
\end{figure}

\begin{figure}[H]
\begin{center}
\includegraphics[height=8cm]{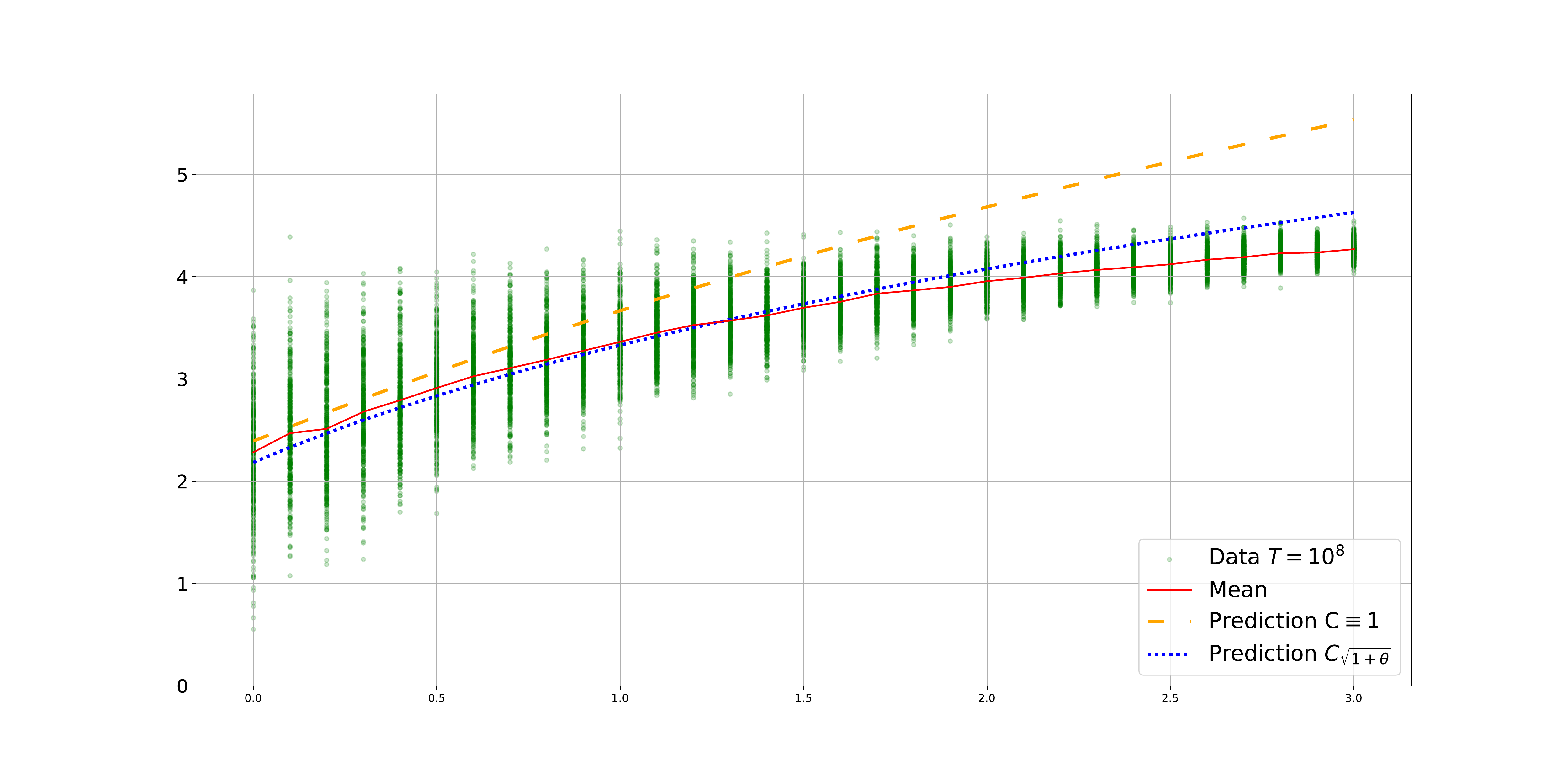}
\end{center}
\caption{The empirical  mean of the samples of $\max_{|h|\leq \pi(\log T)^\theta}\log |\zeta(1/2+\ii(\tau+h))|$ as a function of $\theta$, $0\leq \theta\leq 3$, with step size $0.1$ at $T=10^8$. The dotted lines correspond to the theoretical predictions for $C\equiv 1$ and for $C_{\sqrt{1+\theta}}$.}
\label{fig: mean 8}
\end{figure}

\begin{figure}[H]
\begin{center}
\includegraphics[height=8cm]{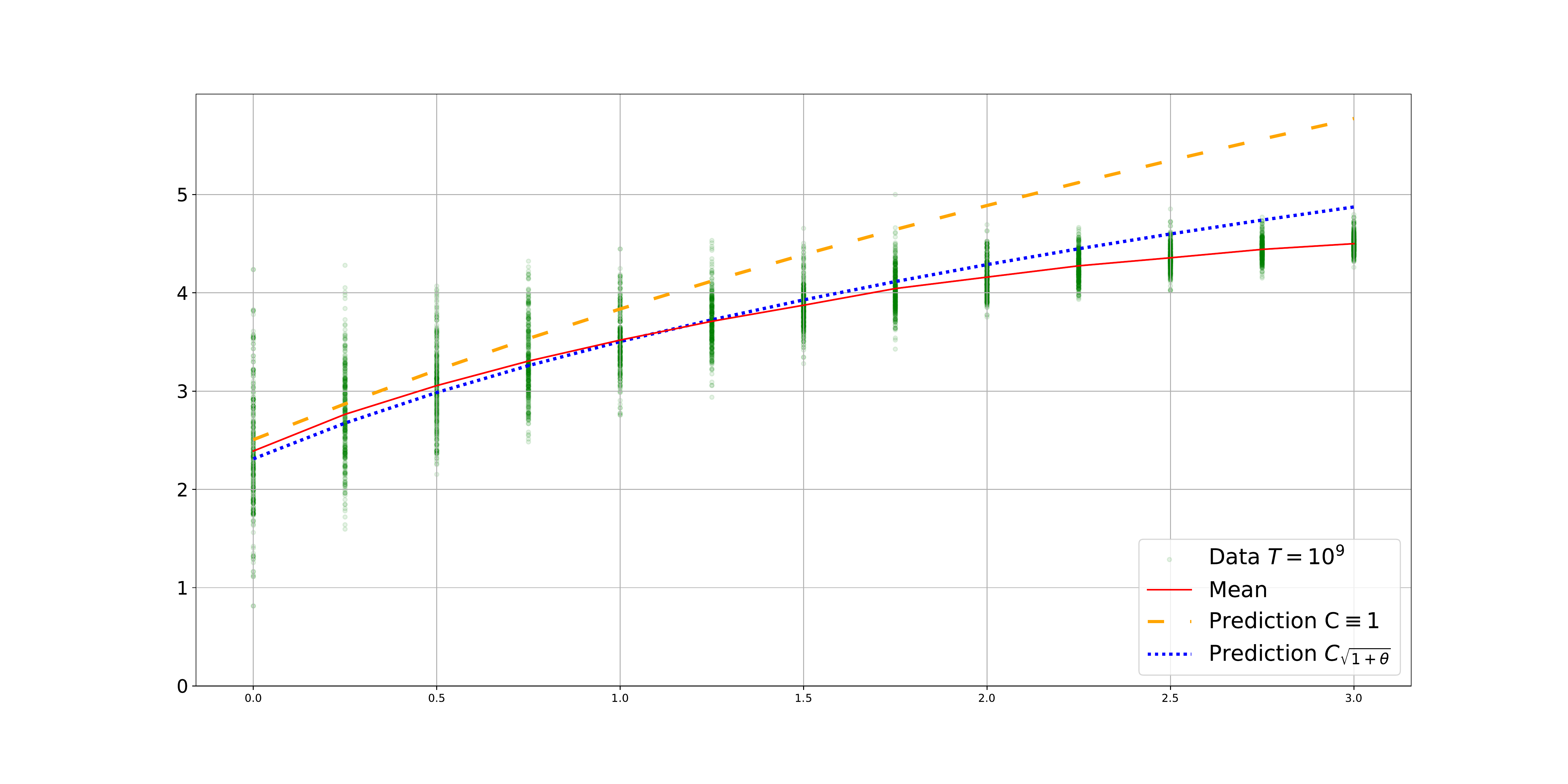}
\end{center}
\caption{The empirical  mean of the samples of $\max_{|h|\leq \pi(\log T)^\theta}\log |\zeta(1/2+\ii(\tau+h))|$ as a function of $\theta$, $0\leq \theta\leq 3$, with step size $0.25$ at $T=10^9$. The dotted lines correspond to the theoretical predictions for $C\equiv 1$ and for $C_{\sqrt{1+\theta}}$.}
\label{fig: mean 9}
\end{figure}

The main numerical results concern the empirical mean of $\max_{|h|\leq \pi(\log T)^{\theta}} \log |\zeta(1/2+\ii(\tau+h))|$ as a function of $0\leq \theta \leq 3$.
These are plotted in Figures \ref{fig: mean 7}, \ref{fig: mean 8} and \ref{fig: mean 9}.
The results are compared with the theoretical predictions of Conjecture \ref{conj: IID} (as detailed in Section \ref{sect: theoretical prediction}), both including the correction $C_{\sqrt{1+\theta}}$ and without, i.e., $C\equiv 1$. We observe that the prediction line for $C\equiv 1$ exhibits greater divergence from the mean as $\theta$ grows, when compared to the corrected $C_{\sqrt{1+\theta}}$ prediction. This is despite the drastic reduction in the variance. This reduction in variance is consistent with the prediction of the parameter $\beta$ of the Gumbel distribution in Equation \eqref{eqn: m beta}, which decreases with $\theta$.
There is a small discrepancy for large $\theta$, where the prediction is slightly outside the range of the sample. The difference between the prediction and the maximum of the sample at $\theta=3$ is  $0.144$ for $T=10^7$ (a relative error of $3\%$), for $T=10^8$ of $0.082$ (a relative error of $1.8\%$), and for $T=10^9$ of $0.076$ (a relative error of $1.6\%$).

Figure \ref{fig: table} gives the ratio of the empirical mean over the two predictions ($C\equiv 1$ and $C_{\sqrt{1+\theta}}$) for $T=10^7$, $T=10^8$ and $T=10^9$.
Again, the ratios suggest that $C_{\sqrt{1+\theta}}$ is the correct prediction. Note that the ratios are improving as $T$ increases.

\begin{figure}[H]
\begin{center}
\includegraphics[height=8cm]{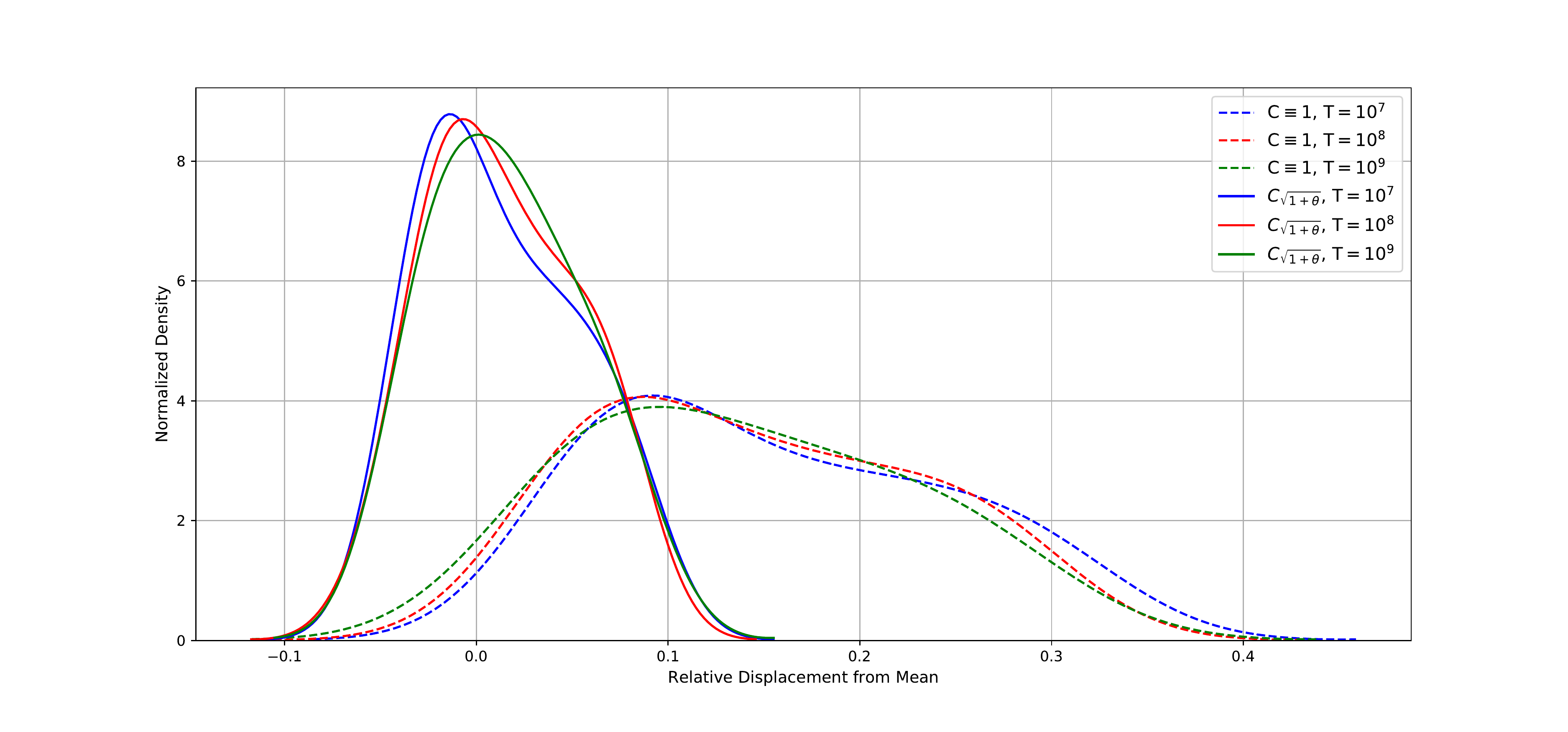}
\end{center}
\caption{The relative displacement of the samples of $\max_{|h|\leq \pi(\log T)^\theta}\log |\zeta(1/2+\ii(\tau+h))|$ from the mean.}
\label{fig: displacement}
\end{figure}

We also examine the convergence of the empirical mean by computing the relative displacement for the predictions from the empirical means, and by calculating the normalized kernel density estimator as shown in Figure \ref{fig: displacement}. We see that as $T$ grows the relative displacements exhibit smaller deviation centered around 0, with the most pronounced effect occurring for $C_{\sqrt{1+\theta}}$ at $T = 10^9$. 

\begin{figure}[H]
\begin{center}
\includegraphics[height=8cm]{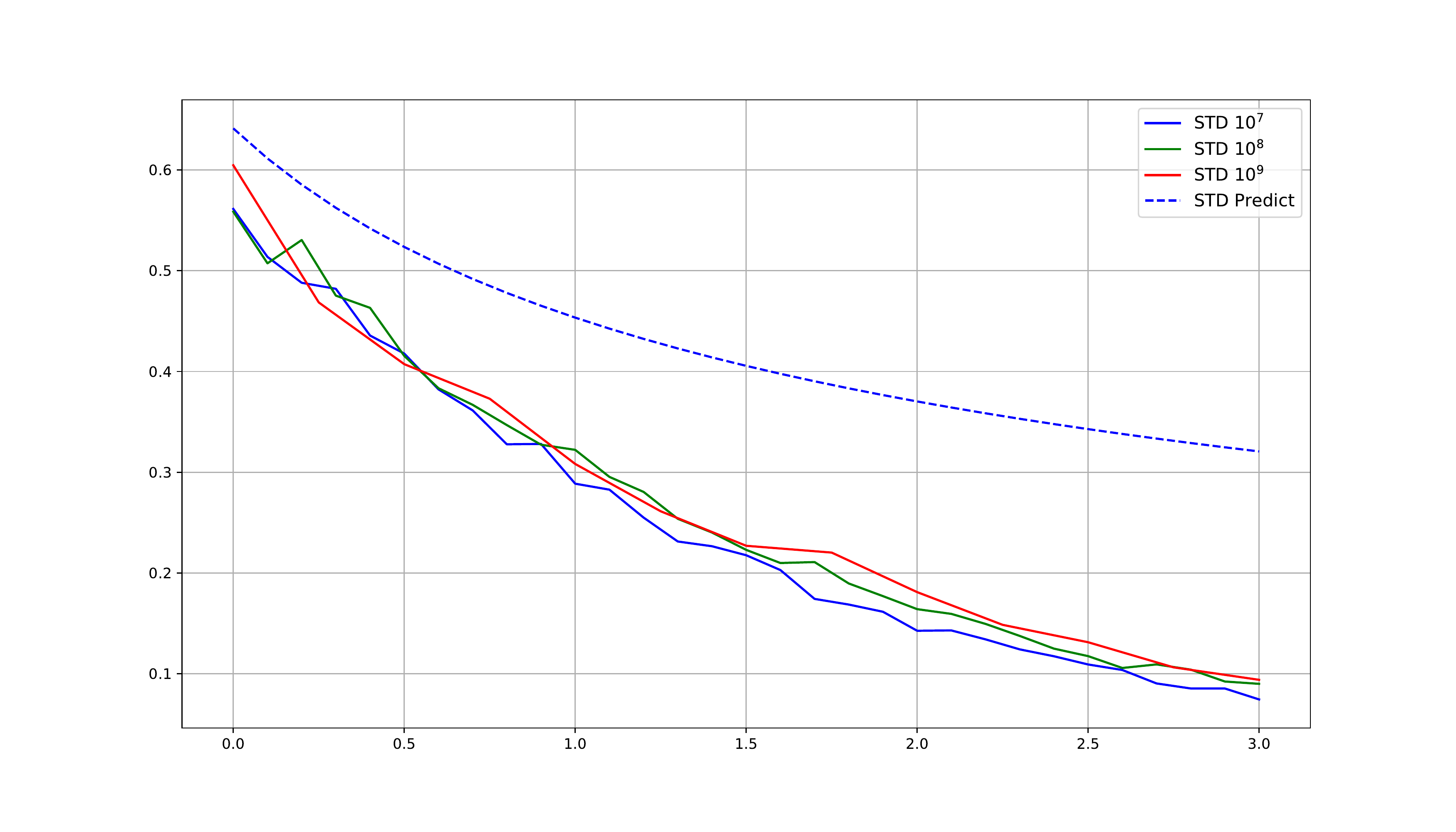}
\end{center}
\caption{The standard deviation of the samples of $\max_{|h|\leq \pi(\log T)^\theta}\log |\zeta(1/2+\ii(\tau+h))|$ as a function of $\theta$, $0\leq \theta\leq 3$. }
\label{fig: std}
\end{figure}

The estimate of the standard deviation of $\max_{|h|\leq \pi(\log T)^{\theta}} \log |\zeta(1/2+\ii(\tau+h))|$ turned out to be trickier, see Figure \ref{fig: std}. The method to obtain the theoretical prediction is explained after Equation \eqref{eqn: std gumbel}. There is a significant discrepancy between the prediction and the numerical results.
We do observe a reduction of the variance as predicted by Equation \eqref{eqn: std gumbel} and the definition of $\beta$ in Equation \eqref{eqn: m beta}.
Note that the standard deviation is fairly small on the whole range of $\theta$. In fact, it is of the order of $1/\sqrt{\mathcal S}$, where $\mathcal S$ is the size of the sample. This might complicate the detection of a signal. The discrepancy seems to be the same for all range of $T$'s. It is also increasing in $\theta$. We currently have no convincing explanations for this phenomenon.

\begin{figure}[H]
\begin{center}

\begin{tabular}{lll}
\begin{tabular}{c| c | c  }
\toprule
$\theta$  & $C\equiv 1$ &   $C_{\sqrt{1+\theta}}$ \\
\midrule
$0$ & 0.9441  & 1.0490 \\
$1$ &  0.9143  &1.0147\\
$2$ & 0.8343 & 0.9679\\
$3$ & 0.7569 & 0.9165\\
\bottomrule
\end{tabular}
\hspace{0.25cm}
&
\hspace{0.25cm}
\begin{tabular}{c| c | c  }
\toprule
$\theta$  & $C\equiv 1$ &   $C_{\sqrt{1+\theta}}$\\
\midrule
$0$ &0.9544 &  1.0454 \\
$1$ &  0.9170  & 1.0099\\
$2$ & 0.8450 & 0.9708\\
$3$ & 0.7713 & 0.9225\\
\bottomrule
\end{tabular}
\hspace{0.25cm}
&
\hspace{0.25cm}
\begin{tabular}{c| c | c  }
\toprule
$\theta$  & $C\equiv 1$ &   $C_{\sqrt{1+\theta}}$\\
\midrule
$0$ & 0.9540  &  1.0343 \\
$1$ &  0.9174  & 1.0043\\
$2$ & 0.8510 & 0.9703\\
$3$ & 0.7795 & 0.9234\\
\bottomrule
\end{tabular}
\end{tabular}

\end{center}
\caption{The ratio of the empirical mean of $\max_{|h|\leq \pi(\log T)^{\theta}} \log |\zeta(1/2+\ii(\tau+h))|$ divided by the model prediction with no correction ($C\equiv1$) and with $C_{\sqrt{1+\theta}}$ at integer $\theta$. From left to right, the data corresponds to $T=10^7$, $T=10^8$, and $T=10^9$.}
\label{fig: table}
\end{figure}

\bibliographystyle{alpha}
\bibliography{zeta_simul_bib}

\end{document}